\documentclass{rspublic}
\usepackage{amsmath}
\usepackage{amsfonts}
\usepackage{amssymb}
\usepackage[Symbol]{upgreek}
\usepackage[nointegrals]{wasysym}
\usepackage{mathrsfs}
\usepackage{graphicx}
\usepackage{accents}
\newcommand{\real}{{\mathbb R}}
\newcommand{\bfu}{{\bf u}}

\newcommand{\bfx}{{\bf x}}
\newcommand{\bfy}{{\bf y}}

\newcommand{\tdel}{\tilde{\delta}}
\newcommand{\Upo}{{\Omega}}
\newcommand{\bdy}{{\partial \Omega}}
\begin{document}

\title[3D Euler Equations]{\large On regularity of the Euler equations in fluid dynamics}

\author[F. Lam]{F. Lam}

\affiliation{ }

\label{firstpage}

\maketitle

\begin{abstract}{Euler Equations; Vorticity Equation}
We assert that the solutions to the Cauchy problem of the inviscid vorticity equation remain regular and unique for any smooth initial data of finite energy. However, the primitive formulation of the Euler equations is not well-posed, due to the passive pressure. One of the implications is that the anomalous energy dissipation, anticipated by Onsager (1949), cannot occur in inviscid flows. In the complete absence of viscous effects, the ultimate accumulation of enstrophy in sealed domains is bound to become arbitrarily excessive, if there is a sustained supply of shears and strains. 
\end{abstract}
\section{Introduction}
In the inviscid thoeyr for incompressible flows, the Euler equations (Euler 1755) are
\begin{equation} \label{euler}
	 \partial_t \bfu  + (\bfu . \nabla) \bfu = - {\rho}^{-1} \nabla p,\;\;\; \nabla.\bfu=0,
\end{equation}
where $\bfu(\bfx,t)=(u,v,w)(\bfx,t)$ is the velocity ($\bfx=(x,y,z)$), the scalar $p(\bfx,t)$ the pressure,  and $\rho$ the density of the fluid. For convenience, we also use the tensor notation, $\bfu(\bfx,t)=u_i(\bfx,t)$, and $\bfx=(x_i)$, $i=1,2,3$.

The scalar pressure can be eliminated from the system \ref{euler}) by making use of the vorticity, ($\nabla {\times} \bfu$), $\upomega(\bfx,t)=(\xi,\eta,\zeta)(\bfx,t)$, which describes the local rotation of fluid elements at any fixed time $t$. Taking curl on (\ref{euler}), we obtain the vorticity equation 
\begin{equation} \label{vort}
	\partial_t \upomega = (\upomega . \nabla) \bfu  - (\bfu . \nabla )\upomega.
\end{equation}
The vorticity field inherits the velocity solenoidal property $\nabla.\upomega=0$. 
As a consistent theoretical framework, the incompressibility hypothesis demands the vorticity-velocity compatibility, 
\begin{equation} \label{uv}
	\Delta \bfu = - \nabla{\times}\upomega.
\end{equation}
In the case of the whole space, the velocity can be found with the aid of the Newtonian potential
\begin{equation} \label{uvt}
	\bfu(\bfx)=\frac{1}{4 \pi}\int_{\real^3} \nabla {\times}\Big(\frac{1}{|\bfx-\bfy|}\Big) \upomega(\bfy) \: \rd \bfy.
\end{equation}
The initial solenoidal data are assumed smooth 
\begin{equation} \label{ic}
 \upomega_0(\bfx)=\nabla{\times}\bfu_0(\bfx)\; (\in C^{\infty}(\Upo)). 
\end{equation} 
In the case of $\real^3$, they have compact supports. We are interested in finite-energy initial value problems subject to the initial data (\ref{ic}). To avoid excessive technicality, we will restrict our discussion to the whole space $\real^3$ or finite domains $\Upo$ with $C^2$-boundary $\bdy$. In inviscid flows, the normal velocity on $\bdy$ remains zero
\begin{equation} \label{nv}
	\bfu(\bfx,t) \cdot \vec{n} (\bfx) = 0, \;\;\; \bfx \in \bdy,
\end{equation}
where $\vec{n}$ denotes the outward normal. 

Taking divergence of (\ref{euler}) and using the continuity, we find that the pressure can be recovered from known shears and strains,
\begin{equation} \label{ppoi}
	\Delta p(\bfx)= - \rho \sum^{3}_{i,j=1} \frac{\partial u_j}{\partial x_i} \: \frac{\partial u_i}{\partial x_j}(\bfx).
\end{equation}
\section{Regularity}

We consider the mollified vorticity field in space and in time:
\begin{equation} \label{mol}
	 \phi_{\tdel}(\upomega) (\bfx,t) = \frac{1}{\tdel^4} \int_{\real^4} \phi \Big( {\bfx'}/{\tdel}, \: {t'}/{\tdel} \Big) \: \upomega (\bfx - \bfx',t - t') \:\rd \bfx' \rd t',
\end{equation}
where the kernel $\phi(\bfx,t) \in C^{\infty}_c $. It follows that the mollified vorticity equation is  
\begin{equation} \label{vortm}
	\partial_t \upomega = (\phi_{\tdel}(\upomega) . \nabla) \bfu  - (\bfu . \nabla )\upomega.
\end{equation}
Adding the three components in (\ref{vortm}), we obtain
\begin{equation} \label{tv}
	\partial_t V = (\phi_{\tdel}(\upomega).\nabla)U - (\bfu.\nabla)V  = R,\;\;\;(\mbox{say})
\end{equation}
where 
\begin{equation*}
	V = \xi + \eta + \zeta, \;\;\; \mbox{and} \;\;\; U = u + v + w,
\end{equation*}
are the total vorticity and velocity respectively. Integrating equation (\ref{tv}) over $\real^3$ or $\Upo$, and making use of the solenoidal property of $\bfu$ and $\upomega$, we establish the vorticity invariance
\begin{equation} \label{ivt}
\frac{\rd }{\rd t} \int_{\Upo} V(\bfx,t) \rd \bfx = 0.
\end{equation}
The initial vorticity, $V(\bfx,0)=(\xi + \eta + \zeta)(\bfx,0)$, is given by data (\ref{ic}). The symmetry or the algebraic cancellation of $R$ at $\bdy$ does not depend on the boundary conditions. In particular, the invariance is independent of the mollification parameter $\tdel$. For any multi-index $\alpha > 0$, we take $\alpha$ derivative of dynamics (\ref{tv}) to get
\begin{equation} \label{tv-alfa}
	{\partial_t (\partial_{x_i}^{\alpha} V)} = \partial_{x_i}^{\alpha} R.
\end{equation}
Because of the space integral over $\Omega$  and the operator $\partial_{x_i}^{\alpha}$commute, the summation of the total vorticity renders the right-hand side ($\partial_{x_i}^{\alpha}{R}$) to zero. Thus we have the extended invariance
\begin{equation} \label{tva}
\frac{\rd }{\rd t} \int_{\Upo} \partial_{x_i}^{\alpha} V(\bfx,t) \rd \bfx = 0.
\end{equation}
Similarly, we observe that 
$	\partial_t (\partial_t^{\beta} V)  = \partial_t^{\beta} R $
for any integer $\beta > 0$. Then the analogous invariance reads
\begin{equation} \label{tvb}
\frac{\rd }{\rd t} \int_{\Upo} \partial_t^{\beta} V(\bfx,t) \rd \bfx = 0.
\end{equation}
In summary, we derive the {\it a priori} bound
\begin{equation} \label{vtlp}
	V(\bfx,t) \in L^p(\bfx \in \Upo)\:L^p(0 \leq t \leq T),\;\;\;1 \leq p < \infty
\end{equation}
for every $T < \infty$. Moreover, the total vorticity is just a finite sum of three real numbers where the components must be essentially space-time bounded
\begin{equation} \label{vteb}
	V(\bfx,t) \in L^{\infty}(\bfx \in \Upo)\:L^{\infty}(0 \leq t \leq T).
\end{equation}
In view of the Sobolev embedding, we conclude that each component of the vorticity is smooth
\begin{equation} \label{vtbd}
	\xi(\bfx,t),\; \eta(\bfx,t),\; \zeta(\bfx,t) \in C^{\infty}(\{\bfx \in \Upo\} \cup \{0 \leq t \leq T \}).
\end{equation}
In $\real^3$, the vorticity now has compact supports. By the standard estimates for the second order elliptic equations, see, for instance, Gilbarg \& Trudinger (1998), we deduce from (\ref{uv}) that
\begin{equation} \label{dulp}
\bfu,\;\nabla \bfu \in C^{\infty},
\end{equation}
which suggests that the incompressibility hypothesis makes sense. In fact, every term, $\partial u_i/\partial x_i$, is bounded and smooth at each instant of time in finite energy flow. Alternatively, the continuity, $\nabla.\bfu=0$, is a weaker restriction on the regularity, as it defines a (finite) zero-sum. The regularity on the pressure follows, once the right-hand side (\ref{ppoi}) is known.

Let $\upomega_1$ and $\upomega_2$ be two vorticity solutions. The corresponding solenoidal velocities are $\bfu_1$ and $\bfu_2$. We denote their differences by $\upvarphi=\upomega_1-\upomega_2$, and $\uppsi=\bfu_1-\bfu_2$. Then the vorticity equation yields
\begin{equation} \label{dm}
	\partial_t \upvarphi = (\upvarphi.\nabla)\bfu_1-(\bfu_1.\nabla)\upvarphi
+(\upomega_2.\nabla)\uppsi-(\uppsi.\nabla)\upomega_2.
\end{equation}
Taking the inner $L^2$ product with $\upvarphi$ and simplifying, we have
\begin{equation*}
\frac{1}{2}\frac{\rd}{\rd t} \|\upvarphi\|^2_{L^2(\real^3)}		\leq \big( \|\nabla\bfu_1\|_{L^{\infty}(\real^3)}+ C \|\upomega_2\|_{L^{\infty}(\real^3)}\big) \|\upvarphi\|^2_{L^2(\real^3)} \leq A_0 \:\|\upvarphi\|^2_{L^2(\real^3)},
\end{equation*}
as $\Delta\uppsi=-\nabla{\times}\upvarphi$, and $\|\nabla\uppsi\|_{L^2}
\leq C \|\upvarphi\|_{L^2}$. By virtue of Gronwall's lemma, we have the bound on the enstrophy
\begin{equation*}
	\|\upvarphi\|^2_{L^2(\real^3)}(\cdot,t) \leq \|\upvarphi(\cdot,0)\|^2_{L^2(\real^3)} \: \exp \bigg(2\int_0^T\!\!A_0(\cdot,s) \:\rd s \bigg),
\end{equation*}
where the exponential term is finite by the bounds (\ref{vtbd}) and (\ref{dulp}). If the initial data $\upvarphi(\bfx,0)$ equals to zero, the uniqueness follows, i.e., $\upomega_1(\bfx,t)=\upomega_2(\bfx,t)$ for all times $t>0$.

In simply-connected domains with $C^2$-boundary, only the normal velocity (see boundary condition (\ref{nv})) vanishes on $\bdy$. Extra analyses are required on the second and fourth terms of (\ref{dm}). We notice that
\begin{equation*}
	\int_{\Upo} \upvarphi.(\bfu_1.\nabla)\upvarphi \: \rd \bfx \leq  \|\bfu_1\|_{L^3(\Upo)} \|\upvarphi\|_{L^6(\Upo)}  \|\nabla\upvarphi\|_{L^2(\Upo)} \leq C_2 \: \|\bfu_1\|_{L^3(\Upo)} \|\nabla \upvarphi\|^2_{L^2(\Upo)},
\end{equation*}
in view of the Sobolev inequality ($C_2>0$). Similarly,
\begin{equation*}
	\int_{\Upo} \upvarphi.(\uppsi.\nabla)\upomega_2 \: \rd \bfx \leq \|\upvarphi\|_{L^6(\Upo)} \|\uppsi\|_{L^2(\Upo)} \|\nabla \upomega_2\|_{L^3(\Upo)} \leq C_4 \: \| \nabla \upvarphi\|^2_{L^2(\Upo)} \|\nabla \upomega_2\|_{L^3(\Upo)}.
\end{equation*}
By the Poincar\'e inequality, the dynamics (\ref{dm}) becomes
\begin{equation*}
\frac{\rd}{\rd t} \|\upvarphi\|^2_{L^2(\Upo)}	\leq A(t)\: \|\upvarphi\|^2_{L^2(\Upo)},
\end{equation*}
where $A = 2\big(\: C'_1\: \|\nabla\bfu_1\|_{L^{\infty}} - C'_2 \: \|\bfu_1\|_{L^3} + C'_3 \: \|\upomega_2\|_{L^{\infty}} - C'_4 \: \|\nabla \upomega_2\|_{L^3}\:\big)< \infty$. Thus the enstrophy is controlled by 
\begin{equation*}
	\|\upvarphi\|^2_{L^2(\Upo)}(\cdot,t) \leq \|\upvarphi(\cdot,0)\|^2_{L^2(\Upo)}\:\exp \bigg(\int_0^T \!\! A(\cdot,s) \: \rd s\bigg),
\end{equation*}
which, in turn, asserts the uniqueness of the vorticity, as well as  the flow-field.

The Euler equations are time-wise reversible and constitute a Hamiltonian system. If there be a blow-up time $t^*$ when the energy approaches an infinite amount, in reverse, we cannot generate a backward motion from $t^*$ where the infinite energy is somehow distributed and sustained. Thus a scenario of finite-time singularity is obscured in physics. For finite initial data in finite domains $\Upo$ where there is no dissipative mechanism, intuition suggests that the integral,
\begin{equation*}
	\int_0^T \big\| \upomega (\bfx,t) \big\|_{L^2(\Upo)} \rd t,
\end{equation*}
increases in time and will go out of bounds for $T \rightarrow \infty$, assuming that there is a continuous supply of shears and strains. Thus, the concept of a {\it global} regularity in inviscid flows must be circumstantial. 
\section{Problems in periodic domains}
Let us consider inviscid flows in a periodic cube: $0 \leq x,y,z \leq L$. The periodic boundary condition is expressed as
\begin{equation*} 
\bfu(\bfx)\;=\;\bfu(\bfx+L).
\end{equation*}
The specification of the velocity allows it to be a strong function of time during the flow evolution. Put simply, any time-dependent velocity field is solenoidal and satisfies the periodic condition. A spurt function, $(\hat{\bfu}, \hat{p})$, is a set of regular solutions such that, for every solution $(\bfu,p)$, the superposition, $(\bfu+\hat{\bfu}, p+\hat{p})$ is also a solution. The following choices work:
\begin{equation*}
\begin{split}
\hat{u}(\bfx,t) & = c_1 \: f(t),\;\;\hat{v}(\bfx,t) = c_2 \: g(t),\;\;\hat{w}(\bfx,t) = c_3 \: h(t),\\
	\hat{p}(\bfx,t)/\rho& = c_1 \: (2L-x)\: f'(t) + c_2 \: (2L-y)\: g'(t) + c_3 \: (2L-z)\: h'(t),
	\end{split}
\end{equation*}
where constants, $c_1,c_2,c_3$, are arbitrary but finite. The spurt functions, $f,g,h$, are smooth, bounded, and arbitrarily chosen so that the energy of the spurts is finite for all times. We impose the initial conditions
\begin{equation*}
	f(0)=f'(0)=0,\;\;\; g(0)=g'(0)=0,\;\;\; \mbox{and} \;\;\;h(0)=h'(0)=0,
\end{equation*}
so that the spurts do not interfere with the initial data $(\bfu,p)(\bfx,0)$. In essence, the Euler equations are ill-posed in the periodic setting. This is an extension of non-uniqueness of the Navier-Stokes equations in periodic domains with periodic boundary conditions. Then both theories are unable to afford well-defined statistical means or ensemble averages. Lastly, it is instructive to mention the fact that there exist no vorticity boundary conditions on $\bdy$. We draw our attention to a regularity criterion of Beale, Kato \& Majda (1984), where an assumption of periodic boundaries has been made.
\section{Non-uniqueness of primitive formulation} 
In the three geometries considered in this section, each of these solutions can be superimposed on to any claimed full solutions of the Euler equations, so that the combined solutions are non-unique.
\subsection{The whole space $\real^3$}
Let a vorticity be
\begin{equation*}
	\xi(\bfx,t)=\gamma(t) \; f_1(y,z),\;\;\;\eta(\bfx,t)=\gamma(t) \; f_2(z,x),\;\;\;\zeta(\bfx,t)=\gamma(t) \; f_3(x,y),
\end{equation*}
where $\gamma$ is an arbitrary regular function, and $\gamma(t=0)=1$, which, effectively, fixes the initial velocity, $\bfu_0$. We assume that smooth and bounded functions, $f_1$ to $f_3$, are compactly supported. Thus, the solenoidal velocity solution, $\bfu(\bfx,t)$, is completely specified by (\ref{uv}) and (\ref{uvt}), at every given time $t$. In spite of the arbitrary velocity field, the Euler equations are indeed satisfied if we equate the pressure gradient $\nabla p$ to
\begin{equation} \label{dp}
	 -\rho \big(\:\partial_t \bfu + (\bfu.\nabla)\bfu \:\big). 
\end{equation}
In other words, the inviscid momenta and energy are conserved at all times. It is instructive to note that the use of the vorticity is a convenient way to define the solenoidal velocity field; the actual vorticity has played no part in the primitive formulation.
\subsection{Flows with fixed tangential velocity in half-space}
In the domain $\Upo_c$: $-\infty < x,y < \infty$, $z\geq0$, and the velocity component $w$ on the $xy$-plane vanishes. To demonstrate the non-uniqueness of the primitive Euler dynamics, we consider
a vorticity function
\begin{equation*}
	\xi(\bfx,t)=\gamma(t)\: \frac{x}{r_c^3},\;\;\;\eta(\bfx,t)=\gamma(t)\:\frac{y}{r_c^3},\;\;\;\zeta(\bfx,t)=\gamma(t)\: \frac{z+c}{r_c^3},
\end{equation*}
where $r_c^2=x^2+y^2+(z+c)^2$, and $c>0$. For finite $a$ and $b$, the velocity boundary conditions are set as
\begin{equation*}
\begin{split}
u(x,y,0,t)& = 2 \pi b^2y\exp\big(-a^2x^2-b^2y^2\:\big) = 2 \pi b^2 u_0(x,y),\\
v(x,y,0,t)& = -2 \pi a^2x\exp\big(-a^2x^2-b^2y^2 \:\big) = - 2 \pi a^2 v_0(x,y),\\
w(x,y,0,t)& = 0,\\
u,v,w& \rightarrow 0,\;\;\;\mbox{as}\;\;\; |\bfx| \rightarrow \infty.
\end{split}
\end{equation*}
The choice of the exponential functions is made without any preference, as the wall tangential components are not fixed in inviscid flows. Extreme boundary data may be prescribed as time-dependent functions, for instance, we may choose $a=a(t)$ or $b=b(t)$. It follows that we have a velocity field
\begin{equation*}
\begin{split}
 u(\bfx,t)&=\int_{\Upo_c} \!\! G_c(\bfx,\bfx') (\partial_{y'} \zeta-\partial_{z'} \eta)(\bfx')  \:\rd \bfx' + b^2z \int_{\real^2} \!\!u_0/R^3 \: \rd x' \rd y',\\
	v(\bfx,t)&=\int_{\Upo_c} \!\! G_c(\bfx,\bfx') (\partial_{z'} \xi-\partial_{x'} \zeta)(\bfx')  \:\rd \bfx'-a^2z \int_{\real^2} \!\!v_0/R^3 \: \rd x' \rd y',\\
	w(\bfx,t)&=\int_{\Upo_c} \!\! G_c(\bfx,\bfx') (\partial_{x'} \eta-\partial_{y'} \xi)(\bfx') \:\rd \bfx',
\end{split}
\end{equation*}
where $G_c(\bfx,\bfx')=(1/R_- - 1/R_+)/(4 \pi)$, $R^2_{\pm}=(x{-}x')^2+(y{-}y')^2+(z {\pm} z')^2$, and $R^2=(x{-}x')^2+(y{-}y')^2+z^2$.
Hence the pair, $(\bfu,\nabla p)$ (cf. (\ref{dp})), are exact and arbitrary solutions of the primitive Euler equations.
\subsection{Multiplicity in a cube}
We give a further example in the box 
$\Upo_b: 0 \leq x \leq L_1, 0 \leq y \leq L_2$, $0 \leq z \leq L_3$. The following regular vorticity distribution is chosen:
\begin{equation*}
	\xi(\bfx,t)=\gamma(t)\: \frac{x+a_1}{r_b^3},\;\;\;\eta(\bfx,t)=\gamma(t)\:\frac{y+a_2}{r_b^3},\;\;\;\zeta(\bfx,t)=\gamma(t)\: \frac{z+a_3}{r_b^3},
\end{equation*}
where $r_b^2=(x+a_1)^2+(y+a_2)^2+(z+a_3)^2$, and $a_1,a_2,a_3>0$. We look for the velocity which vanishes at the walls, $\bfu(\bfx,t)|_{\partial \Upo_b}=0$. By the procedures of the preceding case, the specific velocity is given by
\begin{equation*}
	\bfu(\bfx,t) = \int_{\Upo_b} G_b(\bfx,\bfx') \nabla{\times}\upomega(\bfx') \: \rd \bfx',
\end{equation*}
where the Green's function is expressed as the expansions
\begin{equation*}
\begin{split}
	\frac{8}{\pi^2{\mbox{vol}}(\Upo_b)} \sum_{l,m,n=1}^{\infty} \Big(\: & \sin(l \pi x/L_1)\sin(l \pi x'/L_1) \:\sin(m \pi y/L_2)\sin(m \pi y'/L_2) \\
	\quad & \sin(n \pi z/L_3) \sin(n \pi z'/L_3) \:\Big) \Big/ \Big(\:l^2/L_1^2+m^2/L_2^2+n^2/L_3^2\:\Big).
\end{split}
\end{equation*}
In particular, we set $\gamma(0)=\gamma'(0)=0$ so that the initial data of any full primitive solutions are not affected. 
\section{Discussion}
In fully-developed turbulence in periodic domains, it has been a fallacious view that the Navier-Stokes solutions converge to those of the Euler equations strongly enough in high-Reynolds number flows. In the limit of $Re=UL/\nu \rightarrow \infty$, where $U$ is a representative velocity of a typical large-scale eddy of size $L$, the mean rate of energy dissipation has been assumed to be finite and non-zero. Onsager (1949) conjectured an anomalous dissipation in the complete absence of viscous effects, as long as the velocity field is sufficiently regular, more precisely, if $\bfu \in C^{\kappa}$, where the H\"older exponent $\kappa$ is greater than $1/3$. Probably, the law of energy conservation is violated in irregular velocity fields with $\kappa \leq 1/3$. Our analysis asserts that the anomalous dissipation is impossible because, given smooth bounded initial data, the inviscid vorticity is regular and remains finite at all subsequent times; the derivatives never deteriorate to the postulated rough states.

The dissipation scales in real incompressible flows do not have to reduce indefinitely in the limit $\nu \rightarrow 0$, so that the continuum assumption remains valid, or the small scales are larger than the order of the molecular mean-free path. Briefly, the Reynolds number is {\it not} an appropriate parameter in the description of the dissipation. For example, a flow at $Re = 5{\times} 10^9$ ($\nu=10^{-5} \:\texttt{m$^2$/s}$, $U=25\:\texttt{m/s}$, $L=2000\:\texttt{m}$) must have very different dissipative structures compared with a fluid motion of $U=1\:\texttt{m/s}$ whose linear dimension $L=5\:\texttt{m}$ with a viscosity $\nu = 10^{-9}\:\texttt{m$^2$/s}$. Evidently, we are dealing with two fluids with distinct physical properties, that remain unchanged in constant-temperature flows.
\begin{acknowledgements}
\vspace{2mm}

\noindent 
16 April 2019

\noindent 
\texttt{f.lam11@yahoo.com}

\end{acknowledgements}
%
%\label{lastpage}
\end{document}